\documentclass{iise}
\usepackage{bbm}
\usepackage{amsmath}
\usepackage{multirow}
\usepackage{graphicx}
\usepackage{float}
\usepackage{array}
\usepackage{booktabs}
\usepackage{caption}
\usepackage{subcaption}
\usepackage{cleveref}
\usepackage{bbm}
\usepackage{booktabs}
\usepackage{array}
\usepackage{ragged2e}   
\usepackage{mathtools}  

\usepackage[section]{placeins}
\crefname{equation}{Eq.}{Eqs.}
\crefname{figure}{Fig.}{Figs.}
\crefname{section}{Sec.}{Secs.}
\crefname{table}{Tab.}{Tabs.}

\conference{ } 

\title{\titlesize Empirical Evaluation of Policy-Based Reinforcement Learning for Dynamic Service Control in an M/M/1 Queue}

\author{
Joseph Walton and Gabriel Nicolosi \\ Department of Engineering Management and Systems Engineering\\
Missouri University of Science and Technology, Rolla, MO 65409, USA
}
\authorlist{Walton and Nicolosi}
\abstractID{15941}

\begin{document}
\maketitle

\begin{abstract}

{\small
While reinforcement learning has been increasingly applied to stochastic control, few studies have systematically examined policy-based methods in queuing environments modeled as a semi-Markov decision process (SMDP). To address this gap, we investigate how policy-based reinforcement learning (RL) algorithms perform when applied to the control of service rates in an M/M/1 queue, a common queuing model for manufacturing, computing, and service systems. The problem is formulated as an SMDP in which decisions occur at each new service, allowing an agent to select different service rates from a finite set of speeds, aiming to minimize an objective function that manages system congestion and energy costs. Three policy-based reinforcement learning algorithms, namely REINFORCE, Actor-Critic (A2C), and Proximal Policy Optimization (PPO), are trained in a simulated environment using two state representations: the instantaneous queue length and an augmented state that includes a one-step queue history. Performance is evaluated in terms of convergence speed, sampling efficiency, policy quality, and pseudo-regret relative to the steady-state optimum. 
}

\end{abstract}

\section*{Keywords}
Queuing Control, Semi-Markov Decision Process, Stochastic Control, Sequential Decision Making, Simulation Optimization

\section{Introduction}
Dynamic control of single-server queuing systems is a well-studied problem in operations research. Early studies examined service rate switching policies using queuing-control and dynamic-programming formulations, establishing that optimal policies often exhibit monotone or threshold-type structure under certain objective function and stability assumptions \cite{c69,s70,gh01,tc72}. These results are naturally placed within the Markov decision process (MDP) or semi-Markov decision process (SMDP) frameworks, which allow for average-reward objective functions and event-driven decision epochs \cite{p94}. While such models admit analytically tractable solutions when system parameters are known, uncertainty or increased system complexity can make the use of classical methods impractical, motivating the use of model-free reinforcement learning approaches \cite{zjp24,dg22}. Still, little empirical work has systematically compared standard policy-gradient methods on an $M/M/1$ service rate control problem (under a semi-Markov decision process formulation) with respect to steady-state performance and learning efficiency.

This work is concerned with an empirical study of three average-reward analogues to RL algorithms: average-reward policy gradient (differential REINFORCE) \cite{w92}, actor-critic (A2C) \cite{kt99}, and proximal policy optimization (PPO) \cite{swdrk17}, applied to a dynamic, discrete service rate-controlled $M/M/1$ queue. Performance of each algorithm is then assessed in terms of convergence speed, sampling efficiency, policy quality, and queue-length pseudo-regret relative to a tuned (optimal) threshold-based controller (originally proposed in \cite{tc72}). 

The remainder of this paper is organized as follows: Section~\ref{sec: Related Literature} visits some of the relevant literature on dynamic control of single-server queues and learning-based service rate control. Section~\ref{sec: Problem Construction} presents the problem formulation, defining the controlled $M/M/1$ system as a semi-Markov decision process with an average-reward objective function. Section~\ref{sec: Methodology} describes the reinforcement learning algorithms, experimental design, and evaluation metrics. Section~\ref{sec: Computational Results} presents empirical results and comparative analyses. Finally, Section~\ref{sec: Conclusion and Future Work} concludes with a discussion of the findings and directions for future research.

\section{Related Literature}\label{sec: Related Literature}
Dynamic service rate control for single-server queues has been studied extensively in operations research. Early work showed that, under standard cost structures, optimal policies for queues with controllable service rates are monotone and often take a threshold form \cite{c69, s70, tc72, gh01}. These results were later formalized within an MDP framework supporting average-reward objective functions, which provided a foundation for queuing control analysis \cite{p94}. 

When system parameters are unknown or the system dynamics become analytically intractable, reinforcement learning becomes a practical means of discovering optimal policies. Early policy-gradient methods, namely REINFORCE, demonstrated that parametrized policies can be optimized directly from stochastic system trajectories \cite{w92}. Incorporating value-function baselines leads to actor-critic algorithms, which can improve learning efficiency while retaining on-policy updates suitable for steady-state control at the expense of adding bias \cite{kt99}. More recent methods, such as proximal policy optimization, further improve learning stability through constrained policy updates \cite{swdrk17}.

Recent work has applied reinforcement learning to increasingly complex, high-dimensional queuing networks. Notably, deep reinforcement learning has been applied to multi-server control, achieving strong empirical performance in scenarios where classical methods are computationally prohibitive \cite{dg22}. Extensions of service rate control models have also incorporated reliability effects, such as server breakdowns, while retaining interpretable control structures \cite{zjp24}. Parallel research has introduced regret-based performance metrics for controlled queues, inspired by classical multi-arm bandit problems, linking queue-length regret to learning efficiency and long-run performance guarantees \cite{bc12,ssm21}.

\section{Problem Construction}\label{sec: Problem Construction}
We consider a controlled M/M/1 queue, with arrival and service rates $\lambda$ and $\mu_i$, respectively. The system controller (agent) can select rates from a finite discrete set
\(
\mathcal{A}=\{\mu_1,\mu_2, \dots,\mu_M \},
\)
where $\mu_i$ represents the $i^{th}$ service rate, with $i=1$ representing the smallest available rate (slowest service speed). We assume that $\mu_1 \le \mu_2 \le \dots \le \mu_M$. To ensure system stability in all feasible control actions, the arrival rate satisfies $\lambda < \mu_1$.

Two distinct state representations are used: the first considers only the instantaneous queue length, $s_k = (Q_k)$,
while the second augments this with a second feature of a one-step history of the previous queue length, \textit{i.e.}, $s_k = (Q_k, Q_{k-1})$, attempting to reflect the agent's perception on queue evolution. These two variants allow for comparison between a minimal Markov representation and a partially history-aware formulation.

Decision epochs are enumerated as $k = 1,\dots,K$, where $K$ denotes the number of decision epochs that occur over the time horizon $[0,T]$. These epochs occur immediately before the beginning of each new service. The time epoch each decision occurs at is denoted by $\tau_k \in [0, T]$. Based on a state $s_k$, representative of the queue length, at each decision epoch $k$ the controller selects an action $a_k \in \mathcal{A}$. The chosen service rate remains in effect until the next start of service at the time $\tau_{k+1}$. Between decisions, the system evolves as a $M/M/1$ queue operating at a constant rate $\mu_k$. The time interval between decisions is given by $\Delta \tau_k = \tau_{k+1} - \tau_k$. 

After each decision epoch, the controller receives a (negative) reward reflecting the cost of congestion and energy consumption. The reward at epoch $k+1$ is expressed as a function of the state and chosen service rate at epoch $k$,

\begin{equation}
r_{k+1}=r(s_k,a_k)=- \int_{\tau_k}^{\tau_{k+1}} \big[c_qQ(t)+c_e\mu_k\mathbbm{1}_{\text{busy}} \big] dt,
\end{equation}

where $c_q$ is the congestion cost per unit of time and $c_e$ is an energy scaling factor penalizing high service rates. The objective of the controller is to identify the policy $\pi$ that maximizes the expected long-term average reward:

\begin{equation}
\rho(\pi^*) = \rho^* = \max_{\pi}\lim_{T \to \infty} \frac{1}{T} \mathbb{E}_{\pi}\!\left[
\sum_{k=1}^{K-1} r(s_k, a_k) \right].
\end{equation}

\section{Methodology}\label{sec: Methodology}
Three policy-based algorithms are implemented and trained to learn service policies that minimize congestion and energy costs. The following subsections detail the simulation environment, learning algorithms, baseline policy used for comparison, and performance metrics.

\subsection{Simulation Environment}
All algorithms are evaluated in the controlled $M/M/1$ queuing environment described in Section 3, implemented as a discrete-event simulation where reward is accumulated over each nondeterministic interval between decisions and provided to the agent together with the elapsed time, enabling average-reward learning in the underlying SMDP. To reduce bias from initializing episodes at an empty queue, the initial state of each episode is sampled uniformly from an empirical distribution of the final $n$ states observed in the preceding episodes.

The system parameters were defined by a fixed arrival rate $\lambda=0.0400$, a discrete set of available service rates $\mu_i \in \{0.0417,0.0500,0.0625,0.0833,0.1000\}$, a holding cost coefficient $c_q=0.4$, and an energy cost coefficient $c_e=0.25$. To ensure experimental reproducibility and statistical significance, all experiments were run across five independent trials using a fixed set of pseudo-random seeds $(41,72,99,81,52)$ in \textsf{NumPy}. 

\subsection{Learning Algorithms}
This study evaluates the average-reward (differential) formulations of three policy-based reinforcement learning algorithms: REINFORCE, Actor-Critic (A2C), and Proximal Policy Optimization (PPO). Each method learns a (parameterized) stochastic policy $\pi_{\theta}(a_k|s_k)$, that maps the observed state $s_k$ to a distribution over available service rates $a_k \in \mathcal{A}$ with the goal of improving the long-term (steady-state) performance of the system rather than short-horizon returns. This is accomplished through iteratively updating the policy parameters $\theta$ of a neural network.

Differential REINFORCE, or average-reward policy gradient, updates the policy directly using sampled system trajectories. Following the policy gradient theorem for average rewards, the update rule is defined as:

\begin{equation}
    \theta_{k+1} \doteq \theta_k + \alpha\,(G_k - \rho)\nabla_\theta \ln \pi_\theta(a_k \mid s_k),
\end{equation}

where $\alpha$ is the step size and $G_k$ represents the realized differential return. In this formulation, the return is centered by the current estimate of the long-run average reward $\rho$, resulting in a differential return. Actions that result in rewards exceeding the current estimate of the long-run average reward are reinforced, while actions with below-average outcomes are discouraged \cite{w92}. Although this approach is unbiased, it relies entirely on full-trajectory returns and therefore exhibits high variance, particularly in stochastic queuing environments with variable decision intervals.

A method that mitigates the high variance of differential REINFORCE is A2C, which introduces a learned value function $V_{\phi}(s)$ that estimates the relative desirability of each state \cite{kt99}. The critic $V_\phi$ is parameterized as a two-layer fully connected neural network with hyperbolic tangent activations and a linear output node. The policy parameters are updated using the temporal difference (TD) error $\delta_k$ instead of the realized differential return $G_k$:

\begin{equation} 
    \theta_{k+1} \doteq \theta_k + \alpha_\theta \delta_k \nabla_\theta \ln \pi_\theta(a_k \mid s_k). 
\end{equation}

In the average-reward setting, this error is calculated as:

\begin{equation}
    \delta_k = R_{k+1} - \rho + V_\phi(s_{k+1}) - V_\phi(s_k),
\end{equation}

Here, $\delta_k$ compares the immediate reward $R_{k+1}$ and the change in state value against the average reward $\rho$. By substituting the trajectory-based $G_k$ with the localized error estimate, A2C achieves higher sample efficiency, though it becomes susceptible to approximation errors within the critic $V_{\phi}$.

PPO builds upon the actor-critic framework by explicitly constraining the magnitude of policy updates to promote stability in learning. This is achieved by maximizing a clipped surrogate objective:

\begin{equation}
    L^{\text{CLIP}}(\theta) = \mathbb{E}_k \left[ \min \big( r_k(\theta) \widehat{A}_k,\ \operatorname{clip}(r_k(\theta), 1-\epsilon, 1+\epsilon) \widehat{A}_k \big) \right].
\end{equation}

The term $r_k(\theta) = \pi_{\theta}(a_k|s_k) / \pi_{\theta_{\text{old}}}(a_k|s_k)$ represents the likelihood ratio between the new and old policies. By clipping this ratio within a range defined by $\epsilon>0$, PPO prevents overly aggressive parameter updates that could destabilize the learning process \cite{swdrk17}. The policy parameters are updated by performing multiple epochs of minibatch stochastic gradient ascent on $L^{\text{CLIP}}(\theta)$ using trajectories collected from the current policy. In this study, $\widehat{A}_k$ is computed using the differential average to maintain consistency with the average-reward objective function, providing stability in the presence of the $M/M/1$ queue's inherent variability.

\subsection{Baseline: Threshold-Based Control Policy}
Learned policies are compared against a deterministic threshold-based control scheme defined as

\begin{equation}
\pi^* =
\begin{cases}
\mu_1, & Q < \omega_1,\\[4pt]
\mu_j, & \omega_{j-1} \le Q < \omega_j,\quad j = 2,\dots,M-1,\\[4pt]
\mu_M, & Q \ge \omega_{M-1},
\end{cases}
\end{equation}

where $\pi^*$ denotes the steady-state optimal threshold policy obtained via dynamic programming. The threshold structure is derived from classical results for service rate control in an $M/M/1$ queue, which shows that under standard monotonicity and stability conditions, the optimal policy is non-randomized and monotone in the queue length, admitting one service rate per contiguous interval of queue lengths \cite{tc72}.
The specific threshold values $\omega_i$ are obtained by computing an optimal stationary policy $\pi^*$ via average-reward relative value iteration on a truncated state space $Q\in\{0,\dots,Q_{\max}\}$ \cite{p94}. The resulting dynamic programming policy exhibits monotone switching in $Q$, where switching points are recorded as thresholds, and the final threshold policy is verified by simulation under the same cost function and decision-epoch conventions as the learning environment.

\subsection{Performance Evaluation Metrics}
To assess the empirical performance of each algorithm, four primary metrics are computed: convergence speed, sampling efficiency, policy quality, and queue-length pseudo-regret. Each metric is evaluated over multiple random seeds, and results are reported as the mean and standard deviation across seeds.

\begin{table}[H]
\centering
\caption{Performance metrics used for evaluating learning and control policies.}
\label{tab:metrics}
\begin{tabular}{p{0.18\linewidth} p{0.40\linewidth} p{0.347\linewidth}}
\toprule
\textbf{Metric} & \textbf{Definition} & \textbf{Equation} \\
\midrule

\text{Convergence Speed} ($U_{\eta}$)
&
Minimum number of gradient updates required for the moving average return to reach a target performance threshold $\eta$. &
$U_{\eta} = \min \left\{ U :
\rho_{\theta_{U}}
\ge \rho_\eta \right\}$
\\
\addlinespace[0.3em]

\text{Sampling Efficiency} ($N_\eta$)
&
Minimum number of decision epochs required for the moving-average return to reach a target performance threshold $\eta$.
&
$N_\eta = \min \left\{ N :
\frac{1}{N}\sum_{k=1}^{N}\mathbb{E}[r(s_k,a_k)]
\ge \rho_\eta \right\}$
\\
\addlinespace[0.3em]

\text{Policy Quality}\qquad ($Q_\pi$)
&
Mean reward-rate obtained by the final policy over $K$ evaluation episodes.
&
$Q_\pi=
\frac{1}{K} \sum_{k=1}^{K} \rho(\pi_{k})$
\\
\addlinespace[0.25em]

Queue-Length Pseudo-Regret ($\mathrm{R}_Q(N_{\eta})$)
& Cumulative excess queue length relative to the steady-state optimal threshold policy over training.
& 
$\mathrm{R}_Q(N_{\eta})=\sum_{k=1}^{N_{\eta}} \left( \mathbb{E}[Q_k \mid \pi_k]-\mathbb{E}[Q_k \mid \pi^*] \right)$
\\
\bottomrule
\end{tabular}
\end{table}

\section{Computational Results}\label{sec: Computational Results}
We provide the numerical results from the performance of the differential RL algorithms across two state representations: $s_k = (Q_k)$ and $s_k = (Q_k,Q_{k-1})$. Performance is evaluated in terms of convergence speed ($U_\eta$), sampling efficiency ($N_\eta$), and steady-state policy quality ($Q_\pi$), as defined in Table~\ref{tab:metrics}. Aggregate results across random seeds are summarized in Table~\ref{tab:combined_results}.

\begin{table}[ht]
\centering
\caption{Learning efficiency, convergence speed, and steady-state performance under both state representations.}
\label{tab:combined_results}
\small
\setlength{\tabcolsep}{5pt}
\begin{tabular}{lccc ccc ccc}
\toprule
& \multicolumn{3}{c}{\textbf{Convergence Speed} ($U_{\eta}$)}
& \multicolumn{3}{c}{\textbf{Sampling Efficiency} ($N_{\eta}$)}
& \multicolumn{3}{c}{\textbf{Policy Quality} ($Q_{\pi}$)} \\
\cmidrule(lr){2-4} \cmidrule(lr){5-7} \cmidrule(lr){8-10}
\textbf{Algorithm}
& $(Q_k)$ & $(Q_k,Q_{k-1})$ & $\Delta$\%
& $(Q_k)$ & $(Q_k,Q_{k-1})$ & $\Delta$\%
& $(Q_k)$ & $(Q_k,Q_{k-1})$ & $\Delta$\% \\
\midrule
REINFORCE 
& $5.76\times10^{4}$ & $4.18\times10^{4}$ & $-27.4\%$
& $2.53\times10^{7}$ & $1.67\times10^{7}$ & $-33.9\%$
& $-6.327$ & $-6.317$ & $+0.2\%$ \\

A2C
& $2.29\times10^{6}$ & $4.91\times10^{6}$ & $+114.4\%$
& $2.29\times10^{6}$ & $4.91\times10^{6}$ & $+114.4\%$
& $-6.326$ & $-6.257$ & $+1.1\%$ \\

PPO
& $2.20\times10^{4}$ & $2.05\times10^{4}$ & $-6.8\%$
& $6.53\times10^{6}$ & $7.09\times10^{6}$ & $+8.6\%$
& $-6.329$ & $-6.318$ & $+0.2\%$ \\
\bottomrule
\end{tabular}
\end{table}

\begin{figure}[ht]
    \centering
    \includegraphics[width=\columnwidth]{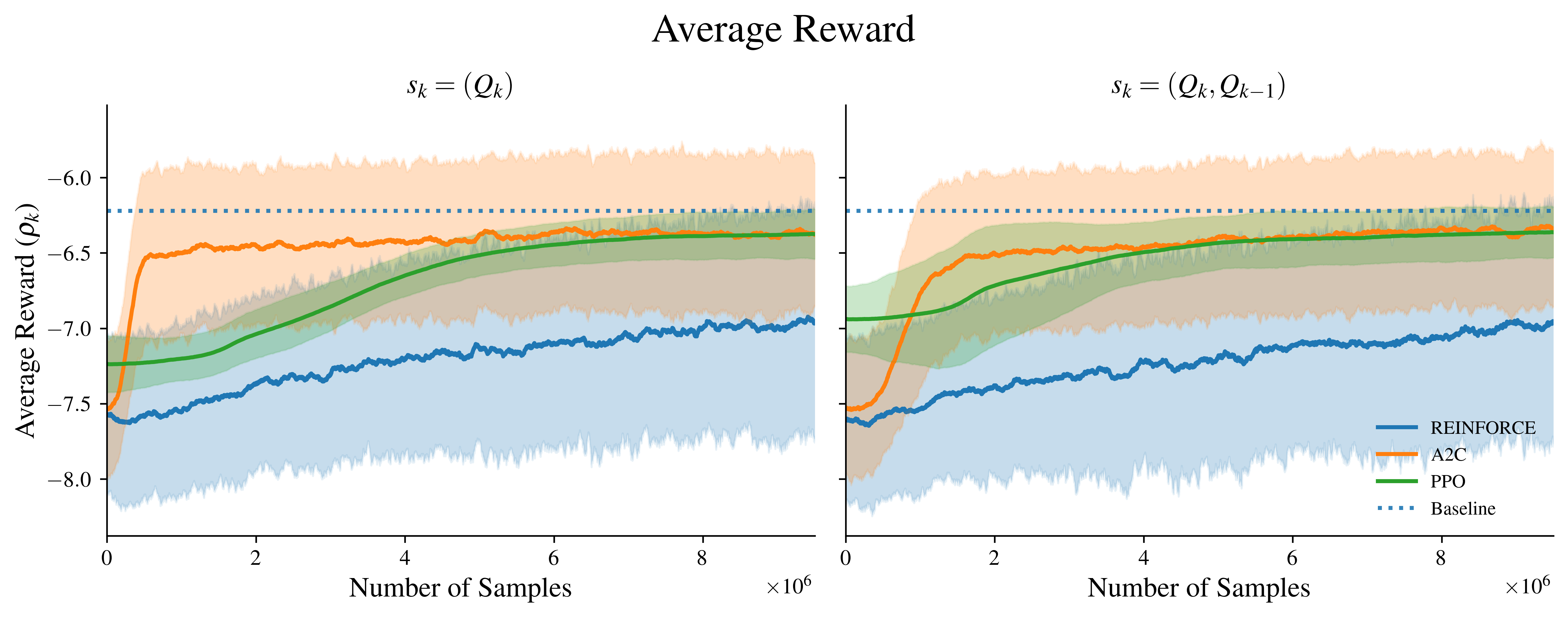}
    \caption{Average reward learning curves versus number of samples for two state representations. Shaded regions indicate $\pm1$ standard deviation across seeds.}
    \label{fig:learning_curve}
\end{figure}

The results from Table~\ref{tab:combined_results} highlight the significant difference in how algorithms handle increased state dimensionality. We note that small changes in early learning dynamics can shift $N_\eta$ even when final performance is similar. PPO emerges as the most computationally efficient method, although the augmented one-step history increased its sampling requirement ($N_{\eta}$) by 8.6\%, its convergence speed ($U_{\eta}$) had improved by 6.8\%. This suggests that PPO’s clipped surrogate objective effectively leverages the richer state information of the one-step history without being destabilized by the noise of the extra feature. Conversely, A2C required more than twice the number of samples (and updates) to reach the performance threshold $\eta$. This change in A2C's learning behavior indicates that critic error propagates more aggressively as dimensionality increases. Interestingly, REINFORCE proved surprisingly receptive to the augmented state, showing a 33.9\% reduction in $N_{\eta}$, while still being the most sample-intensive method by an order of magnitude to reach the performance threshold. The impact of the learning stability difference can be observed in Figure \ref{fig:learning_curve}.

The learning dynamics in Figure \ref{fig:learning_curve} reveal that while all methods converge to a near-optimal steady-state reward rate ($\approx-6.32$), their trajectories differ substantially. PPO exhibited the lowest variance across seeds, whereas A2C showed rapid initial improvement in conjunction with sustained oscillation, followed by REINFORCE, which improved slowly due to its sparse updates.

\begin{figure}[ht]
    \centering
    \includegraphics[width=\columnwidth]{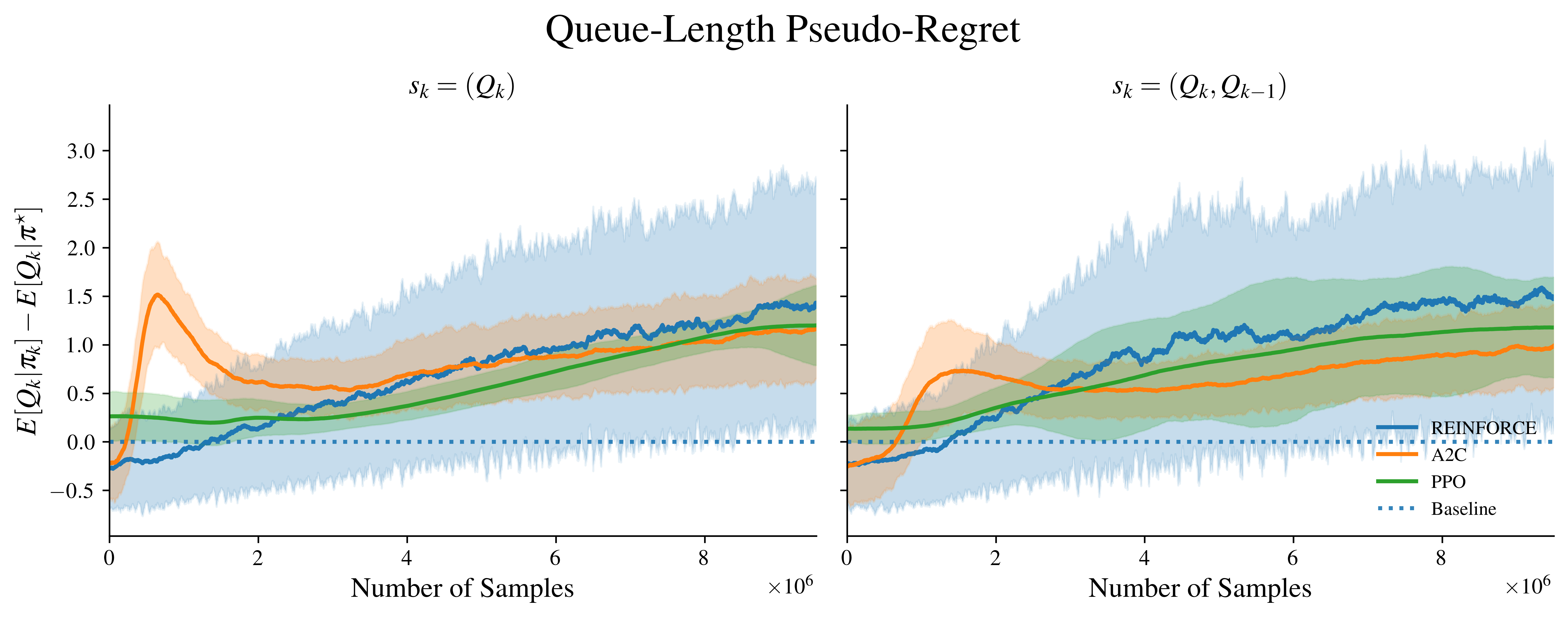}
    \caption{Difference between the expected queue length under the policy at $k$ and the expected queue length under the optimal policy $\pm1$ standard deviation across seeds.}
    \label{fig:queue_len}
\end{figure}

The impact of these efficiency gaps is captured by the Queue-Length Pseudo-Regret in Figure \ref{fig:queue_len}. This metric represents the cumulative excess congestion incurred relative to the optimal threshold policy \cite{ssm21}. PPO maintains the shallowest queue-length regret curve, whereas A2C exhibits a distinct spike in regret in early training. This spike reflects the congestion cost of sub-optimal exploration before the critic stabilizes.

\section{Conclusion and Future Work}\label{sec: Conclusion and Future Work}
The empirical results presented in the paper demonstrate trade-offs between stability, sample efficiency, and computational intensity in service rate control of an $M/M/1$ queue. PPO provides the most reliable balance between sampling efficiency and optimization stability. The primary discovery is the fragility of A2C in the presence of a one-step state history, where a single additional feature doubled the amount of samples required until convergence ($N_{\eta}$). This confirms that for stochastic service systems, policy-clipping mechanisms are of great importance when moving beyond single-feature Markovian states. 

Future work will investigate advanced feature engineering to mitigate noisy state representations of history while providing the agent with richer perceptual information. Key initiatives include showing the difference between states (delta encoding) $\Delta Q_k$ to explicitly express queue dynamics, larger history windows consisting of ($Q_k, Q_{k-1},\dots, Q_{k-n}$) to provide larger trend signals, log-normalized queue-lengths $s_k=\ln(1+Q_k)$ to help the agent distinguish low-buffer states where the risk of idleness is high, and extending the environment to include stochastic server breakdowns. This requires the agent’s perception to include binary availability, testing the algorithm's ability to learn "risk-aware" service policies that pre-clear congestion when the probability of failure increases \cite{zjp24}.

\bibliography{References}

@article{dg22,
	title = {Queueing {Network} {Controls} via {Deep} {Reinforcement} {Learning}},
	volume = {12},
	doi = {10.1287/stsy.2021.0081},
	number = {1},
	journal = {Stochastic Systems},
	author = {Dai, J. G. and Gluzman, Mark},
	year = {2022},
	pages = {30--67},
}

@article{zjp24,
	title = {Adaptive service rate control of an {M}/{M}/1 queue with server breakdowns},
	volume = {106},
	doi = {10.1007/s11134-023-09900-z},
	language = {en},
	number = {1-2},
	journal = {Queueing Systems},
	author = {Zheng, Yi and Julaiti, Juxihong and Pang, Guodong},
	year = {2024},
	pages = {159--191},
}

@article{gh01,
	title = {Dynamic {Control} of a {Queue} with {Adjustable} {Service} {Rate}},
	volume = {49},
	number = {5},
	journal = {Operations Research},
	author = {George, Jennifer M. and Harrison, J. Michael},
	year = {2001},
	pages = {720--731},
}

@article{ssm21,
	title = {Learning {Algorithms} for {Minimizing} {Queue} {Length} {Regret}},
	volume = {67},
	doi = {10.1109/TIT.2021.3054854},
	number = {3},
	journal = {IEEE Transactions on Information Theory},
	author = {Stahlbuhk, Thomas and Shrader, Brooke and Modiano, Eytan},
	year = {2021},
	pages = {1759--1781},
}

@misc{swdrk17,
	title = {Proximal {Policy} {Optimization} {Algorithms}},
	doi = {10.48550/arXiv.1707.06347},
	publisher = {arXiv},
	author = {Schulman, John and Wolski, Filip and Dhariwal, Prafulla and Radford, Alec and Klimov, Oleg},
	year = {2017},
}

@article{bc12,
	title = {Regret {Analysis} of {Stochastic} and {Nonstochastic} {Multi}-armed {Bandit} {Problems}},
	volume = {5},
	number = {1},
	journal = {Foundations and Trends® in Machine Learning},
	author = {Bubeck, Sébastien and Cesa-Bianchi, Nicolò},
	year = {2012},
	publisher = {Now Publishers, Inc.},
	pages = {1--122},
}

@article{w92,
	title = {Simple statistical gradient-following algorithms for connectionist reinforcement learning},
	volume = {8},
	doi = {10.1007/BF00992696},
	number = {3},
	journal = {Machine Learning},
	author = {Williams, Ronald J.},
	year = {1992},
	pages = {229--256},
}

@inproceedings{kt99,
	title = {Actor-{Critic} {Algorithms}},
	author = {Konda, Vijay R and Tsitsiklis, John N},
    booktitle = {Advances in Neural Information Processing Systems},
	year = {1999},
}

@techreport{s70,
    author = {Sabeti, Houshang},
    title = {{Optimal} {Decision} {in} {Queueing}},
    institution = {California University, Berkeley, Operations Research Center},
    year = {1970},
    doi = {0.21236/AD0708007}
}

@phdthesis{c69,
    author = {Crabill, Thomas Buskirk},
    title = {{Optimal} {Control} {of} {a} {Queue} {with} {variable} {Service} {Rates}},
    school = {Cornell University},
    year = {1969}
}

@book{p94,
	series = {Wiley {Series} in {Probability} and {Statistics}},
	title = {Markov {Decision} {Processes}: {Discrete} {Stochastic} {Dynamic} {Programming}},
	shorttitle = {Markov {Decision} {Processes}},
	publisher = {John Wiley \& Sons, Inc},
	author = {Puterman, Martin L.},
	year = {1994},
}

@article{tc72,
	title = {Optimal {Control} of a {Service} {Facility} with {Variable} {Exponential} {Service} {Times} and {Constant} {Arrival} {Rate}},
	volume = {18},
	doi = {10.1287/mnsc.18.9.560},
	number = {9},
	journal = {Management Science},
	author = {Crabill, Thomas B.},
	year = {1972},
	pages = {560--566},}

\end{document}